\documentclass[a4paper,12pt]{article}
\usepackage[T2A]{fontenc}
\usepackage[cp1251]{inputenc}
\usepackage[english]{babel}
\usepackage{amsthm}
\usepackage[tbtags]{amsmath}
\usepackage{amsfonts,amssymb}
\begin{document}

\newtheorem{theorem}{Theorem}
\newtheorem{lemma}{Lemma}
\newtheorem{proposition}{Proposition}
\newtheorem{Cor}{Corollary}

\begin{center}
{\bf Incidence Coalgebras: Automorphisms and Derivations}
\end{center}
\begin{center}
Piotr Krylov\footnote{National Research Tomsk
State University, e-mail: krylov@math.tsu.ru .},
Askar Tuganbaev\footnote{National Research University MPEI; Lomonosov Moscow State University; e-mail: tuganbaev@gmail.com .}
\end{center}

\textbf{Abstract.} We describe automorphisms and derivations of the incidence coalgebra $\text{Co}(X,F)$ of the partially ordered set $X$ over a field $F$. In this case, the fact is significantly used that the dual algebra of the coalgebra $\text{Co}(X,F)$ is isomorphic to the incidence algebra $I(X,F)$ of the partially ordered set $X$ over the field $F$.


\textbf{MSC2020 database 16T15, 16W25}

The study is supported by grants of Russian Science Foundation (=RSF), no.~23-21-00375, https://rscf.ru/en/project/23-21-00375/ (P.A. Krylov) and no.~22-11-00052, https://rscf.ru/en/project/22-11-00052 (A.A. Tuganbaev).

\tableofcontents

\section{Introduction}\label{section1} 

In this paper, for a partially ordered set $X$ and a field $F$, we study automorphisms and derivations of the incidence coalgebra $\text{Co}(X,F)$ of $X$ over $F$. It is known that the algebra dual to the coalgebra $\text{Co}(X,F)$ is canonically isomorphic to the incidence algebra $I(X,F)$ of the partially ordered set $X$ over a field $F$. The theory of incidence algebras and coalgebras is presented in the book \cite{SpiO97}.\label{7} For more information, we can refer to papers \cite{AquS00}\label{1} and \cite{SmiW94} \label{6} which also contain a variety of material on arbitrary coalgebras.
 
The structure of the automorphism group and the derivation space of the incidence algebra $I(X,F)$ has been clarified quite a long time ago. The corresponding results are presented in \cite{SpiO97}.\label{7} Subsequently, these results were transferred and generalized in a fairly large number of papers to more general incidence rings $I(X,R)$ , where $X$ is a pre-ordered set and $R$ is some ring (e.g., see \cite{Khr10},\label{2} \cite{Khr12},\label{3} \cite{KryT23}\label{4}). In a number of articles, rings that are close in different senses to ordinary incidence rings were introduced and studied (e.g., see \cite{LerS81}\label{5} and \cite{SmiW94}\label{6}).

Under the study of automorphisms and derivations, we pass from the coalgebra $\text{Co}(X,F)$ to the algebra $I(X,F)$, and conversely. However, the reverse transition is difficult. In the general case, it is impossible to return from the algebra $I(X,F)$ to the  coalgebra $\text{Co}(X,F)$ in any standard way. At the same time, it is important that automorphisms and derivations of the coalgebra $\text{Co}(X,F)$ induce automorphisms and derivations of the dual algebra and the algebra $I(X,F)$, respectively.

We denote the automorphism group of the partially ordered set $X$ by $\text{Aut }X$. For any two elements $x,y\in X$, we denote by $[x,y]$ the subset $\{z\in X\,|\,x\le z\le y\}$. It is called an \textsf{interval} in $X$. Further, we assume that all intervals in $X$ are finite. In such a case, $X$ is called a \textsf{locally finite} partially ordered set.

If $V$ is a linear space over a field $F$ (an $F$-space, for brevity), then $V^*$ is the dual space of $V$, i.e. $V^*=\text{Hom}_F(V,F)$. Further, $\text{End}_FV$ denotes the endomorphism ring (i.e. the ring of linear operators) of the space $V$ and $\text{Aut}_FV$ denotes the automorphism group (i.e. the group of invertible linear operators) of the space $V$. The role of the space $V$ is played either some algebra $A$ or the coalgebra $C$ and its dual algebra $C^*$.

Let $A$ be an algebra over a field $F$ (briefly, an $F$-algebra). Then $\text{Aut }A$ is the automorphism group of the algebra $A$, $\text{In(Aut }A)$ is the group of inner automorphisms, $\text{Der }A$ is the derivation space of the algebra $A$ and $\text{In(Der }A)$ is a subspace of inner derivations. The automorphism group of the coalgebra $C$ is also denoted by $\text{Aut }C$.

For groups $G$ and $H$, we denote the semidirect product of $G$ and $H$ by $G\leftthreetimes H$.\\

\section*{\textbf{Chapter I.\\ SUPPORTING RESULTS}}\label{chapterI}
\addtocontents{toc}{\textbf{Chapter I. SUPPORTING RESULTS}\par}
\section{Coalgebras and Dual Algebras}\label{section2} 

Let $(C,\Delta,\varepsilon)$ be some coalgebra over a field $F$. Thus, $C$ is an $F$-space, $\Delta$ is a comultiplication in $C$, $\varepsilon$ is a counit for $C$. We use the letter $C$ to denote this coalgebra.

We can naturally define the linear mappings $m\colon C^*\otimes C^*\to C^*$ and $u\colon F\to C^*$ in such a way that we get the resulting $F$-algebra $C^*= (C^*,m,u)$.
It is called the \textsf{dual algebra} for the coalgebra $C$. The identity element of the algebra $C^*$ is a counit $\varepsilon$. We remark that if $C$ is a finite-dimensional space, then the duality also exists in the opposite direction.

We elaborate a little on what has just been said. The multiplication $m\colon C^*\otimes C^*\to C^*$ is induced by the comultiplication $\Delta$. Namely, let $\rho\colon C^*\otimes C^*\to (C\otimes C)^*$ be the canonical monomorphism, where
$$
(\rho(\chi\otimes\xi))(c\otimes d)=\chi(c)\xi(d)
$$
for all $\chi,\xi\in C^*$ and $c,d\in C$. We set $m=\Delta^*\rho$, where $\Delta^*\colon (C\otimes C)^*\to C^*$ is the induced mapping and $\Delta^*(\theta)=\theta\Delta$ for $\theta\in(C\otimes C)^*$. Further, we denote by $\omega$ the canonical isomorphism $F\otimes F\to F$. Then $\rho(\chi\otimes\xi)=\omega\cdot(\chi\otimes\xi)$ for $\chi,\xi\in C^*$.

Let $c\in C$ and $\Delta(c)=\sum_ia_i\otimes b_i$, $a_i,b_i\in C$. We have the relations
$$
(m(\chi\otimes\xi))(c)=(\Delta^*\rho(\chi\otimes\xi))(c)=
(\rho(\chi\otimes\xi)\Delta)(c)=(\omega\cdot(\chi\otimes\xi)\Delta)(c)=
$$
$$
=(\omega\cdot(\chi\otimes\xi))(\sum_ia_i\otimes b_i)=
\omega(\sum_i\chi(a_i)\otimes \xi(b_i))=\sum_i\chi(a_i)\xi(b_i).
$$

We denote by $\chi\circ\xi$ the element $m(\chi\otimes\xi)$. In the following lemma, the multiplication rule in the algebra $C^*$ is given.

\textbf{Lemma 2.1.} Let $\chi,\xi\in C^*$. Further, let $c\in C$ and $\Delta(c)=\sum_ia_i\otimes b_i$, where $a_i,b_i\in C$. Then we have $$
(\chi\circ\xi)(c)=\sum_i\chi(a_i)\xi(b_i).
$$

We denote by $\Gamma$ the mapping $\text{End}_FC\to\text{End}_FC^*$, where $\Gamma(\varphi)=\varphi^*$ for every $\varphi\in \text{End}_FC$. Here $\varphi^*$ is a linear mapping of the space $C^*$ induced by $\varphi$, i.e. $\varphi^*(\chi)=\chi\varphi$ for every $\chi\in C^*$. If $\varphi\in\text{Aut}_FC$, then it is clear that $\varphi^*\in \text{Aut}_FC^*$.

We write down the following useful fact.

\textbf{Lemma 2.2.} The mapping $\Gamma$ is an antimonomorphism of $F$-algebras $\text{End}_FC\to\text{End}_FC^*$ and a group antimonomorphism $\text{Aut}_FC\to\text{Aut}_FC^*$.

Everywhere else $X$ is some locally finite partially ordered set. Let $C$ be a linear space such that its basis consists of all intervals $[x,y]$ of the set $X$. We define mappings $\Delta\colon C\to C\otimes C$ and $\varepsilon\colon C\to F$ by setting
$$
\Delta([x,y])=\sum_{x\le z\le y}[x,z]\otimes [z,y],
$$
$$
\varepsilon([x,y])=
\begin{cases}
1, \text{if } x=y;\\
0, \text{if } x\ne y
\end{cases}
$$
for every interval $[x,y]$.

The triple $(C,\Delta,\varepsilon)$ is a coalgebra which is called the \textsf{incidence coalgebra} of the locally finite partially ordered set $X$. We denote it by $\text{Co}(X,F)$ or, like before, by the letter $C$.

\section{Algebra Dual to Coalgebra $\text{Co}(X,F)$ and Incidence Algebra $I(X,F)$}\label{section3} 

Sometimes we denote an incidence algebra $I(X,F)$ by the letter $A$. The incidence coalgebra $\text{Co}(X,F)$ is denoted by $C$, as agreed at the end of the previous section. We reveal the connections between the algebras $C^*$ and $A$.

There are mutually inverse isomorphisms of $F$-algebras
$$
\Phi\colon C^*\to A \text{ and } \Psi\colon A\to C^*,
$$
$$
\text{where }\, (\Phi(\chi))(x,y) =\chi([x,y]),\; \chi\in C^*,
$$
$$
(\Psi(f))([x,y]) =f(x,y),\; f\in A,
$$
for all $x,y\in X$ with $x\le y$.

In its turn, $\Phi$ and $\Psi$ induce mutually inverse isomorphisms of algebras of linear mappings
$$
\Phi^*\colon \text{End}_FC^*\to \text{End}_FA,
$$
$$
\Psi^*\colon \text{End}_FA\to \text{End}_FC^*,
$$
$$
\text{where }\, \Phi^*(\eta)=\Phi\eta\Psi,\; \eta\in\text{End}_FC^*,
$$
$$
\Psi^*(\xi)=\Psi\xi\Phi,\; \xi\in\text{End}_FA.
$$

We denote by $\Theta$ the antimonomorphism $\Psi^*\Gamma\colon \text{End}_FC\to \text{End}_FA$ of algebras (see Lemma 2.2). Let's find out how $\Theta$ works; this is very useful. For arbitrary $\varphi\in\text{End}_FC$ and $f\in A$, we can write down the relations
$$
(\Theta(\varphi))(f)=((\Phi^*\Gamma)(\varphi))(f)=
(\Phi^*(\Gamma(\varphi)))(f)=((\Phi\Gamma(\varphi))\Psi)(f)=
$$
$$
=\Phi(\Gamma(\varphi)(\Psi(f)))=\Phi(\Psi(f)\varphi), \text{ where } \Psi(f)\varphi\in C^*.
$$

For any two elements $x,y\in X$ with $x\le y$, we have the relation
$$
\Phi(\Psi(f)\varphi)(x,y)=\Psi(f)(\varphi([x,y])).
$$

As a result, we obtain the relation
$$
((\Theta(\varphi))(x,y)=\Psi(f)(\varphi([x,y])).
$$

Further, if
$$
\varphi([x,y])=\alpha_1[x_1,y_1]+\ldots+\alpha_k[x_k,y_k],\, \text{ where } \alpha_i\in F,
$$
then the relations
$$
\Psi(f)(\varphi([x,y]))=\alpha_1\Psi(f)([x_1,y_1])+\ldots+\alpha_k\Psi(f)([x_ky_k])=
$$
$$
=\alpha_1f(x_1,y_1)+\ldots+\alpha_kf(x_k,y_k)
$$
hold.

Thus,, we can write down the following proposition.

\textbf{Proposition 3.1.} Let $\varphi\in\text{End}_FC$, $f\in I(X,F)$, $[x,y]$ is an interval in $X$. There are assertions below.
\begin{enumerate}
\item[\textbf{1.}]
The relation
$$
((\Theta(\varphi))(f))(x,y)=\Psi(f)(\varphi([x,y])).
$$
is true. In particular, we obtain the relation
$$
\Psi(f)([x,y])=f(x,y)\quad \text{(for } \varphi=1).
$$
\item[\textbf{2.}]
If $\varphi([x,y])=\alpha_1[x_1,y_1]+\ldots+\alpha_k[x_k,y_k]$, $\alpha_i\in F$, then the relation
$$
((\Theta(\varphi))(f))(x,y)=\alpha_1f(x_1,y_1)+\ldots+\alpha_kf(x_k,y_k)
$$
holds.
\end{enumerate}

We also obtain the following relations related to multiplication in the dual algebra $C^*$. Of course, they are consistent with Lemma 2.1. Namely, for arbitrary functions $f,g\in A$ and an interval $[x,y]$, we have
$$
(\Psi(f)\circ \Psi(g))([x,y])=(\Psi(fg))([x,y])=fg([x,y])= \sum_{x\le z\le y}f(x,z)g(z,y).
$$

\section{On Group $\text{Aut}(I(X,F))$}\label{section4} 

We put together some information about the incidence algebra $I(X,F)$ and its automorphism group (see \cite{KryT23}\label{4} and \cite{SpiO97}\label{7} for details). Similar to the previous section, the algebra $I(X,F)$ is also denoted by the letter $A$. For the automorphism group of the partially ordered set $X$, we use the symbol $\text{Aut }X$.

We recall some special functions in $I(X,F)$. For a given element $x\in X$, we set $e_x(x,x)=1$ and $e_x(z,y)=0$ for all remaining pairs $(z,y)$. The system $\{e_x\,|\,x\in X\}$ consists of pairwise orthogonal idempotents in the ring $L_1$ (the last ring is defined in the next paragraph).

We define a subring $L_1$ and an ideal $M_1$ in $A$. We set
$$
\begin{cases}
L_1=\{f\in A\,|\,f(x,y)=0,\;\text{if } x\ne y\},\\
M_1=\{f\in A\,|\,f(x,y)=0,\;\text{if } x=y\}.
\end{cases}
$$

We have a direct sum $A=L_1\oplus M_1$ of $F$-spaces. Therefore, the algebra $A$ is a splitting extension of the ideal $M_1$ with the use of the subring $L_1$. The ideal $M_1$ can be considered as an $L_1$-$L_1$-bimodule and a non-unital algebra, as well.

Let we have an arbitrary element $x\in X$. We denote by $R_x$ the set of all functions $f\in A$ such that $f(z,y)=0$ for $(z,y)\ne (x,x)$. The relations $R_x=e_xAe_x=e_xL_1e_x$ hold. Consequently, $R_x$ is a ring with identity element $e_x$. More precisely, $R_x\cong F$.

Now we take two distinct elements $x,y$ and we set
$$
M_{xy}=\{f\in A\,|\,f(s,t)=0,\, \text{ if }\, (s,t)\ne (x,y)\}.
$$

Here $M_{xy}=e_xAe_y$ and, therefore, $M_{xy}$ is an $R_x$-$R_y$-bimodule. The relation $M_{xy}=e_xM_1e_y$ is also true; in addition, $M_{xy}\cong F$.

The product $\prod_{x,y\in X}M_{xy}$ has a natural structure of an $L_1$-$L_1$-bimodule. In this bimodule, we also can define a multiplication by the use of the relation
$$
(g_{xy})(h_{xy})=(d_{xy}),\quad \text{where }\, d_{xy}=\sum_{x\le z\le y}g_{xz}h_{zy}.
$$ 

Thereafter, the product $\prod_{x,y\in X}M_{xy}$ becomes a (non-unital) algebra.

\textbf{Proposition 4.1 \cite[Proposition 3.1]{KryT23}.} There are canonical algebra isomorphisms $L_1\cong\prod_{x\in X}R_{x}$, $M_1\cong\prod_{x,y\in X}M_{xy}$ and a canonical $L_1$-$L_1$-bimodule isomorphism $M_1\cong\prod_{x,y\in X}M_{xy}$.

Let $\varphi$ be an arbitrary automorphism of the algebra $A$. Based on the direct sum $A=L_1\oplus M_1$ of $F$-spaces, we can associate with the automorphism $\varphi$~ $2\times 2$ the matrix
$\left(\begin{smallmatrix}\alpha&\gamma\\ \delta&\beta\end{smallmatrix}\right)$. In addition, $\gamma=0$ in our case (see \cite[Proposition 3.1]{KryT23})\label{4}. Further, $\alpha$ is an automorphism of the algebra $L_1$ and $\beta$ is an automorphism of the (non-unital) algebra $M_1$. If $\delta=0$, then $\beta$ also is an automorphism of an $L_1$-$L_1$-bimodule $M_1$.

Let $v$ be some invertible element of the algebra $A$. With the use of it, we can obtain an automorphism $\mu_v$ of this algebra if we set $\mu_v(f)=v^{-1}fv$ for every $f\in A$. Such an automorphism is called the \textsf{inner automorphism determined by the element $v$}. All inner automorphisms form a normal subgroup $\text{In(Aut }A)$ in $\text{Aut }A$. In our case, it is useful to detail the  defined notion as done in \cite{KryT23}\label{4}. $\text{In}_1(\text{Aut }A)$ (resp., $\text{In}_0(\text{Aut }A)$) are denoted by subgroup of inner automorphisms determined by invertible elements of the form $1+g$, $g\in M_1$ (resp., determined by invertible elements of the algebra $L_1$). The first subgroup is normal in $\text{Aut }A$ and we have a semidirect decomposition
$$
\text{In(Aut }A)=\text{In}_1(\text{Aut }A)\leftthreetimes \text{In}_0(\text{Aut }A).
$$

In addition to $\text{In(Aut }A)$, we define two more subgroups of the group $\text{Aut }A$.

We denote by $\text{Mult }A$ the subgroup of $\text{Aut }A$ consisting of the automorphisms of the form $\left(\begin{smallmatrix}1&0\\ 0&\beta\end{smallmatrix}\right)$. Such automorphisms are said to be \textsf{multiplicative}. Let's approach their definition a little differently.

\textbf{Definition 4.2.} Let we have a system of non-zero elements 
$\{c_{xy}\in F\,|\,x<y\}$ satisfying the relation
$$
c_{xy}=c_{xz}c_{zy}\; \text{ such that }\; x<z<y.\eqno (1)
$$

The systems of elements defined in Definition 4.2 are called \textsf{multiplicative systems}. It was possible to define a multiplicative system as a system $\{c_{xy}\,|\,x\le y\}$, where $c_{xx}=1$ for all $x$ and $c_{xy}=c_{xz}c_{zy}$ such that $x\le z\le y$.

Let $\psi\in\left(\begin{smallmatrix}1&0\\ 0&\beta\end{smallmatrix}\right)\in\text{ Mult }A$. For any two elements $x,y$ with $x<y$, there exists a non-zero element $c_{xy}$ of the field $F$ such that the relation $\beta(g)=c_{xy}g$ is true for all $g\in M_{xy}$. In addition, the system $\{c_{xy}\,|\,x<y\}$ satisfies to the relation $(1)$ (see the paragraph before \cite[Proposition 10.2]{KryT23})\label{4}. Thus, we can associate the automorphism $\psi$ with the multiplicative system $\{c_{xy}\,|\,x<y\}$. Conversely, every multiplicative system $\{c_{xy}\,|\,x<y\}$ determines a multiplicative automorphism $\psi$. Namely, for an element $g=(g_{xy})\in M_1$, we set $\psi(g)=(c_{xy}g_{xy})$ and $\psi(f)=f$ for $f\in L_1$.

We can write down the following fact \cite[Proposition 10.2]{KryT23}\label{4}.

\textbf{Proposition 4.3.} There is a one-to-one correspondence between multiplicative automorphisms and multiplicative systems of elements.

We consider another type of automorphisms. Let $\tau\in\text{Aut }X$. We define a mapping $\eta_{\tau}\colon A\to A$ by setting
$$
\eta_{\tau}(f)(x,y)=f(\tau(x),\tau(y))\; \text{ for all } f\in A,\; x,y\in X.
$$
Here $\eta_{\tau}$ is an automorphism of the algebra $A$ and the correspondence $\tau\to \eta_{\tau}$ is an antiisomorphic embedding $p\colon \text{Aut }X\to\text{Aut }A$ of groups. We identify $\tau$ with $\eta_{\tau}$ using the mapping $p$. The automorphisms of the form $\eta_{\tau}$ are called \textsf{order} automorphisms. The image of the embedding $p$ is denoted by $\text{Aut}_AX$.

We give the main result about the structure of the group $\text{Aut }A$. It is contained in \cite[Theorem 7.3.6]{SpiO97}.\label{7} Some its generaizations and clarification are obtained in \cite[Corollary 9.4]{KryT23})\label{4}.

\textbf{Theorem 4.4.} There are the following relations of groups
$$
\text{Aut }A=(\text{In(Aut }A)\cdot\text{Mult }A)\leftthreetimes \text{Aut}_AX=
$$
$$
=\text{In}_1\text{(Aut }A)\leftthreetimes\text{Mult }A\leftthreetimes \text{Aut}_AX.
$$

\section*{\textbf{Chapter II. AUTOMORPHISM GROUP\\ OF INCIDENCE COALGEBRAS}}\label{chapterII}
\addtocontents{toc}{\textbf{Chapter II. AUTOMORPHISM GROUP\\ 
\mbox{}\hspace{6mm}OF INCIDENCE COALGEBRAS}\par}

\section{Group $\text{Aut }C$}\label{section5} 

Let $C=(C,\Delta,\varepsilon)$ and $C'=(C',\Delta',\varepsilon')$ be two arbitrary coalgebras over a field $F$. A linear mapping $\varphi\colon C\to C'$ is called a \textsf{homomorphism from the coalgebra} $C$ into the coalgebra $C'$ provided that $(\varphi\otimes \varphi)\Delta=\Delta'\varphi$ and $\varepsilon'\varphi=\varepsilon$.

If $\varphi\colon C\to C$ is a homomorphism of coalgebras and $\varphi$ is a bijection, then $\varphi$ is called an \textsf{automorphism} of the coalgebra $C$. All automorphisms of the coalgebra $C$ form a group with respect to composition. It is denoted by $\text{Aut }C$ and is called the \textsf{automorphism group} of the coalgebra $C$. Certainly, $\text{Aut }C$ is a subgroup in the group $\text{Aut}_FC$.

We recall that $\text{Aut }C^*$ is the automorphism group of the dual algebra $C^*$ (see Section 2).

\textbf{Lemma 5.1.} If $\varphi\in \text{Aut }C$, then $\varphi^*\in \text{Aut }C^*$.

\textbf{Proof.} First of all, we have to verify that the relation $\varphi^*m=m(\varphi^*\otimes\varphi^*)$ holds. Since $\varphi$ is an automorphism of the coalgebra $C$, the relation 
$\Delta\varphi=(\varphi\otimes\varphi)\Delta$ holds. We also have the relations
$$
\varphi^*m=\varphi^*\Delta^*\rho=\Delta^*(\varphi\otimes\varphi)^*\rho.
$$
On the other hand, the relation
$m(\varphi^*\otimes\varphi^*)=\Delta^*\rho(\varphi^*\otimes\varphi^*)$ is true. With use of simple calculations, it is easy to verify that the relation $(\varphi\otimes\varphi)^*\rho=\rho(\varphi^*\otimes\varphi^*)$ is true, which is required.

The identity element of $C^*$ is $\varepsilon$. It follows from the relation $\varepsilon\varphi=\varepsilon$ that $\varphi^*(\varepsilon)=\varepsilon\varphi=\varepsilon$, i.e. $\varphi^*$ maps the identity element to the identity element. It is clear that $\varphi^*$ is a bijection. Consequently, $\varphi^*$ is an automorphism of the algebra $C^*$.~$\square$

\section{Two Subgroups in $\text{Aut }C$}\label{section6} 

Until the end of the chapter, $C$ is an incidence coalgebra $\text{Co}(X,F)$ and $A$ is the incidence algebra $I(X,F)$. For the algebra $A$, we define analogues of matrix units. Let $x,y\in X$ with $x<y$. We denote by $e_{xy}$ a function $X\times X\to R$ such that $e_{xy}(s,t)=1$ if $s=x$, $t=y$, and $e_{xy}(s,t)=0$ for all other pairs $(s,t)$. The functions $e_{xy}$ satisfy the following property:
$$
\text{if } x<z<y, \text{ then } e_{xz}e_{zy}=e_{xy}.
$$

In view of Lemmas 2.2, 5.1 and the material in Section 3, we have group antimonomorphisms
$$
\Gamma\colon \text{Aut }C\to\text{Aut }C^* \text{ and }\; \Theta\colon \text{Aut }C\to\text{Aut }A.
$$
Now we define subgroups in $\text{Aut }C$ which are similar to the subgroups $\text{Mult }A$ and $\text{Aut}_AX$ in $\text{Aut }A$ from Section 4.

Let we have a multiplicative system of elements $\{d_{xy}\,|\,x< y\}$ (similar systems appeared in Section 4). We define a mapping $\lambda\colon C\to C$ by setting $\lambda([x,y])=d_{xy}[x,y]$ for every basis vector $[x,y]$ ($x<y$) of the space $C$, and $\lambda([x,x])=[x,x]$ for every vector $[x,x]$. We can naturally extend $\lambda$ to a linear mapping of the space $C$. It is easy to verify that $\lambda$ is an automorphism of the coalgebra $C$. We call it the \textsf{multiplicative automorphism} corresponding to the multiplicative system $\{d_{xy}\,|\,x< y\}$.

All multiplicative automorphisms of the coalgebra $C$ form a subgroup in $\text{Aut }C$ which is denoted by $\text{Mult }C$.

\textbf{Proposition 6.1.} The groups $\text{Mult }C$ and $\text{Mult }A$ are isomorphic and the isomorphism carried out by the mapping $\Theta$.

\textbf{Proof.} First, we remark that the groups $\text{Mult }A$ and $\text{Mult }C$ are Abelian.

Let $\lambda\in\text{Mult }C$ and let $\{d_{xy}\,|\,x< y\}$ be the corresponding multiplicative system. We have $\lambda([x,y])=d_{xy}[x,y]$ for all basis vectors $[x,y]$ of the space $C$ (we can assume that $d_{xx}=1$ for all vectors of the form $[x,x]$). With the use of Proposition 3.1 for any pair $(s,t)$ with $s\le t$, we can write down the relations
$$
((\Theta(\lambda))(e_{xy}))(s,t)=d_{st}e_{xy}(s,t)=
\begin{cases}
d_{xy},\,\text{ if } (s,t)=(x,y);\\
0,\,\text{ if } (s,t)\ne(x,y).
\end{cases}
$$
Therefore, we obtain $(\Theta(\lambda))(e_{xy})=d_{xy}e_{xy}$. We conclude that $\Theta(\lambda)$ is a multiplicative automorphism of the algebra $A$ corresponding to the system $\{d_{xy}\,|\,x< y\}$ (we have also take into account Proposition 4.3). Therefore, $\Theta(\lambda)\in\text{Mult }A$.

Now we take some automorphism $\psi\in\text{Mult }A$. Let $\psi$ correspond to the multiplicative system $\{c_{xy}\,|\,x< y\}$ from Proposition 4.3. We denote by $\delta$ the multiplicative automorphism of the coalgebra $C$ corresponding to the system $\{c_{xy}\,|\,x< y\}$. It follows from the previous paragraph that $\Theta(\delta)=\psi$. It can be argued that the restriction of $\Theta$ to $\text{Mult }C$ maps from $\text{Mult }C$ onto $\text{Mult }A$ and, therefore, is an isomorphism $\text{Mult }C\to\text{Mult }A$.~$\square$

Now we consider an analogue of order automorphisms for the coalgebra $C$. Let we have an automorphism $\tau$ of a partially ordered set $X$. We define a mapping $\tau\colon C\to C$ (we use the same symbol $\tau$ to designate it) by setting $\tau([x,y])=[\tau(x),\tau(y)]$ for every basis vector $[x,y]$. We obtain a linear mapping $\tau$ of the space $C$ with converse mapping $\tau^{-1}$. More precisely, $\tau$ is an automorphism of the coalgebra $C$.

We also define a group monomorphism $q\colon\text{Aut }X\to\text{Aut }C$, where $q(\tau)=\tau$, $\tau\in \text{Aut }X$. The image of the mapping $q$ is denoted by $\text{Aut}_CX$. After Proposition 4.3,  the group antimonomorphism $p\colon\text{Aut }X\to\text{Aut }A$ appeared.

\textbf{Proposition 6.2.} \textbf{1.} The relation $\Theta q=p$ is true.

\textbf{2.} The mapping $\Theta$ is a group antiisomorphism $\text{Aut}_CX\to \text{Aut}_AX$.

\textbf{Proof.} \textbf{1.} Let $\tau\in \text{Aut }X$, $f\in A$, and let $x,y\in X$ with $x<y$. We clarify that
$$
((p(\tau))(f))(x,y)=f(\tau(x),\tau(y)),\quad
q(\tau)([x,y])=[\tau(x),\tau(y)].
$$
By the use of Proposition 3.1, we now write down the relations
$$
(((\Theta q)(\tau))(f))(x,y)=((\Theta(q(\tau)))(f))(x,y)=
f(\tau(x),\tau(y))=((p(\tau))(f))(x,y).
$$
Therefore, $\Theta q=\rho$.

\textbf{2.} The assertion directly follows from \textbf{1}.~$\square$

With use of Proposition 6.2, we obtain the relation $\Gamma q=\Psi^*p$, as well; see Section 3.

\section{Inner Automorphisms of Coalgebra $\text{Co}(X,F)$}\label{section7} 

We introduce the notion of an inner automorphism of the coalgebra. The purpose of the section is to verify that the group of inner automorphisms of the coalgebra $\text{Co}(X,F)$ is antiisomorphic to the group of inner automorphisms of the algebra $I(X,F)$.

For a while let $C$ be an arbitrary coalgebra.

\textbf{Definition 7.1.} An automorphism $\nu$ of the coalgebra $C$ is said to be \textsf{inner} if $\nu^*$ is an inner automorphism of the dual algebra $C^*$. 

All inner automorphisms of the coalgebra $C$ form a normal subgroup in $\text{Aut }C$ which is denoted by $\text{In(Aut }C)$.

\textbf{From this point until the end of the section, the letter $C$ again denotes the incidence coalgebra $\text{Co}(X,F)$ and $A$ is the incidence algebra $I(X,F)$.}

Let $\nu\in \text{In(Aut }C)$. We have an inclusion $\nu^*\in \text{In(Aut }C^*)$. Then we can to verify that $\Phi^*(\nu^*)\in \text{In(Aut }A)$. Therefore, we obtain a group antiisomorphic embedding
$$
\Theta\colon \text{In(Aut }C)\to \text{In(Aut }A)
$$
(see the beginning of Section 6).

Based on an arbitrary invertible function $h$ of the algebra $A$, we define some elements of the field $F$. We take specific elements $x,y$ of the set $X$ with $x\le y$. For elements $s,t\in X$ with $x\le s\le t\le y$, we set
$$
\alpha_{xy}(s,t)=h^{-1}(x,s)\cdot h(t,y).\eqno (1)
$$

\textbf{Lemma 7.2.} The following assertions are true.
\begin{enumerate}
\item[\textbf{1.}]
For elements $x,s,r,t,y\in X$ with $x\le s\le r\le t\le y$, we have the relation
$$
\alpha_{xy}(s,t)=\alpha_{xr}(s,r)\cdot \alpha_{ry}(r,t).
$$
\item[\textbf{2.}]
Let we have $x,p,q,u,v,y\in X$ with $x\le p\le q<u\le v\le y$. Then 
the sum $\sum_{q\le z\le u}\alpha_{xz}(p,q)\cdot\alpha_{zy}(u,v)$ is equal to zero.
\item[\textbf{3.}]
The following relations are true:
$$
\alpha_{xx}(x,x)=1,\quad \sum_{x\le s\le y}\alpha_{xy}(s,s)=0 \text{ for } x<y.
$$
\end{enumerate}

\textbf{Proof.} \textbf{1.} It follows from relation $(1)$ that
$$
\alpha_{xr}(s,r)\cdot \alpha_{ry}(r,t)=h^{-1}(x,s)h(r,r)h^{-1}(r,r)\cdot h(t,y).
$$
It remains to remark that $h(r,r)=1=h^{-1}(r,r)$, see \cite[Theorem 4.3]{KryT23}\label{4}.

\textbf{2.} We can write down the relations
$$
\sum_{q\le z\le u}\alpha_{xz}(p,q)\cdot \alpha_{zy}(u,v)=
\sum_{q\le z\le u}h^{-1}(x,p)\cdot h(q,z)\cdot h^{-1}(z,u)\cdot h(v,y)=
$$
$$
=h^{-1}(x,p)\cdot h(v,y)\left(\sum_{q\le z\le u}h(q,z)\cdot h^{-1}(z,u)\right)=
$$
$$
=h^{-1}(x,p)\cdot h(v,y)\cdot hh^{-1}(q,u)=h^{-1}(x,p)\cdot h(v,y)\cdot 1_A(q,u)=0
$$
(here $1_A$ is the identity mapping of the algebra $A$).

The relations of assertion \textbf{3} are directly verified.~$\square$

As before, $h$ is some invertible function in $A$. We define a linear mapping $\nu$ of the space $C$ by setting
$$
\nu([x,y])=\sum_{x\le s\le t\le y}\alpha_{xy}(s,t)[s,t]=
\sum_{x\le s\le t\le y}h^{-1}(x,s)[s,t]h(t,y).\eqno (2)
$$
for every basis vector of $V$.
 
We denote by $\mu$ the inner automorphism of the algebra $A$ determined by its invertible element $h$.

\textbf{Proposition 7.3.} The mapping $\nu$ is an inner automorphism of the coalgebra $C$. In addition, the relation $\Theta(\nu)=\mu$ is true.

\textbf{Proof.} We verify the relation $\Delta\nu=(\nu\otimes\nu)\Delta$, where, as before, $\Delta$ is a comultiplication in $C$.

We take an arbitrary interval $[x,y]$. By calculating, we obtain the relations
$$
\Delta\nu([x,y])=\sum_{x\le s\le r\le t\le y}\alpha_{xy}(s,t)([s,r]\otimes[r,t]),\eqno (3)
$$
$$
(\nu\otimes\nu)\Delta([x,y])=\sum_{x\le p\le q\le z\le u\le v\le y}\alpha_{xz}(p,q)\alpha_{zy}(u,v)([p,q]\otimes[u,v]).\eqno (4)
$$

It is necessary to verify the coincidence of the coefficients for the same basis vectors of the space $C\otimes C$ contained in the right sides of relations $(3)$ and $(4)$.

There is one-to-one correspondence between the summands in $(3)$ and the summands in $(4)$ with $q=z=u$. The specific summand $\alpha_{xy}(s,t)([s,r]\otimes[r,t])$ in $(3)$ corresponds to the  summand $\alpha_{xr}(s,r)\alpha_{ry}(r,t)([s,r]\otimes[r,t])$ in $(4)$.
In view of Lemma 7.2, $\alpha_{xy}(s,t)=\alpha_{xr}(s,r)\alpha_{ry}(r,t)$.

Now we prove that the sum of all summands in $(4)$ with $q<u$ is equal to zero. To do this, we divide this sum into the sum of summands. Then it will be easy to see that all such summands are equal to zero. 

We take some basis vector $[p,q]\otimes[u,v]$ with $q<u$. There are several summands with this basis vector, since $z$ takes values in the interval $[q,u]$. The coefficient of the taken vector is equal to $\sum_{q\le z\le u}\alpha_{xz}(p,q)\alpha_{zy}(u,v)$. The last expression is equal to zero, by Lemma 7.2. The relation $\Delta\nu=(\nu\otimes\nu)\Delta$ is proved.

We also have to verify the relation $\varepsilon\nu=\varepsilon$. Let $[x,y]$ be an arbitrary basis vector of the coalgebra $C$. If $x=y$, then $\varepsilon([x,x])=1$ and $\nu([x,x])=\alpha_{xx}(x,x)[x,x]$, where $\alpha_{xx}(x,x)=1$, by Lemma 7.2.

If $x<y$, then $\varepsilon([x,y])=0$. By Lemma 7.2, we have $\sum_{x\le s\le y}\alpha_{xy}(s,s)=0$. With the use relation $(2)$, we obtain the required result.

It remains to verify that $\nu$ is a bijection. First, we verify that the relation $\Theta(\nu)=\mu$ is true. Then it will be easy to obtain that $\nu$ is bijective.

Let $f\in A$ and $x,y\in X$ ($x\le y$). With the use of Proposition 3.1 we have the relation
$$
((\Theta(\nu))(f))(x,y)=\sum_{x\le s\le t\le y}\alpha_{xy}(s,t)f(s,t).
$$
We also write down the relation $(\mu(f))(x,y)=h^{-1}fh(x,y)$ and do some calculations:
$$
h^{-1}fh(x,y)=\sum_{x\le t\le y}(h^{-1}f)(x,t)\cdot h(t,y)=
\sum_{x\le t\le y}\left(\sum_{x\le s\le t}h^{-1}(x,s)f(s,t)\right)\cdot h(t,y)=
$$
$$
=\sum_{x\le s\le t\le y}h^{-1}(x,s)\cdot h(t,y)\cdot f(s,t)=
\sum_{x\le s\le t\le y}\alpha_{xy}(s,t)f(s,t).
$$

The relation $\Theta(\nu)=\mu$ is really true.

Similar to $(1)$, we set
$$
\beta_{xy}(s,t)=h^{-1}(t,y)h(x,s).\eqno(5)
$$
Then we define linear mapping $\varkappa\colon C\to C$ by the use of the relation
$$
\varkappa([x,y])=\sum_{x\le s\le t\le y}\beta_{xy}(s,t)[s,t].\eqno(6)
$$
Similar to the case of $\nu$, we can verify that $\varkappa$ is a homomorphism of coalgebras and $\Theta(\varkappa)=\mu^{-1}$. Now it follows from relations $\Theta(\nu\varkappa)=1=\Theta(\varkappa\nu)$ that $\nu\varkappa=1=\varkappa\nu$ and $\nu$, $\varkappa$ are mutually inverse isomorphisms. The last relations are also verified by direct calculation.

Thus, $\nu$ is an inner automorphism of the coalgebra $C$ and $\Theta(\nu)=\mu$. The proof of the proposition is completed.~$\square$

We conclude that all inner automorphisms of the coalgebra $C$ can be obtained from invertible elements of the algebra $A$ with the use of relations $(1)$ and $(2)$. The following statement gives the exact meaning to the last phrase.

\textbf{Corollary 7.4.} The restriction of the mapping $\Theta$ to $\text{In(Aut }C)$ is an antiisomorphism between the groups $\text{In(Aut }C)$ and $\text{In(Aut }A)$.

We set 
$$
\text{In}_0(\text{Aut }C)=\Theta^{-1}(\text{In}_0(\text{Aut }A)),\quad \text{In}_1(\text{Aut }C)=\Theta^{-1}(\text{In}_1(\text{Aut }A)).
$$

Before Definition 4.2, it was write down the semidirect decomposition
$$
\text{In(Aut }A)=\text{In}_1(\text{Aut }A)\leftthreetimes \text{In}_0(\text{Aut }A).
$$

Taking into account Corollary 7.4, we write down the following fact.

\textbf{Corollary 7.5.} We have a semidirect decomposition of groups
$$
\text{In(Aut }C)=\text{In}_1(\text{Aut }C)\leftthreetimes \text{In}_0(\text{Aut }C).
$$

\section{Structure of Group $\text{Aut }C$}\label{section8} 

We have Theorem 4.4 about the structure of the group $\text{Aut }(I(X,F))$ and Proposition 6.1, 6.2 and Corollaries 7.4, 7.5, as well. 

They allow to formulate the main result of Chapter 1.

We remark that the designation $\text{Aut}_CX$ appeared before Proposition 6.2.

\textbf{Theorem 8.1.} Let $A$ be an incidence algebra $I(X,F)$ and let $C$ be the incidence coalgebra $\text{Co}(X,F)$. The following assertions are true.
\begin{enumerate}
\item[\textbf{1.}]
There are the relations
$$
\text{Aut }C=(\text{In}(\text{Aut }C)\text{Mult }C)\leftthreetimes 
\text{Aut}_CX=
\text{In}_1(\text{Aut }C)\leftthreetimes\text{Mult }C \leftthreetimes \text{Aut}_CX.
$$
\item[\textbf{2.}]
The automorphism groups $\text{Aut }C$ and $\text{Aut }A$ are antiisomorphic. For example, the antiisomorphism can be obtained with the use of the mapping $\Theta$.
\end{enumerate}

\textbf{Proof. \textbf{1.}} We take an arbitrary automorphism $\varphi\in\text{Aut }C
$. Then $\Theta(\varphi)\in\text{Aut }A$ (see Lemma 5.1 and the beginning of Section 6). By Theorem 4.4, we have the relation $\Theta(\varphi)=\mu\psi\tau$, where $\mu\in\text{In}_1(\text{Aut }A)$, $\psi\in\text{Mult }A$, $\tau\in \text{Aut}_AX$. Now we use Propositions 6.1, 6.2 and Corollary 7.4 to write down
$$
\Theta(\nu)=\mu,\; \Theta(\lambda)=\psi,\; \Theta(\sigma)=\tau,
$$
where $\nu\in \text{In}_1(\text{Aut }C)$, $\lambda\in\text{Mult }C$, $\sigma\in\text{Aut}_CX$. Therefore, $\Theta(\sigma\lambda\nu)=\Theta(\nu)\Theta(\lambda)\Theta(\sigma)$. Consequently, $\varphi=\sigma\lambda\nu$ or $\varphi=\nu'\lambda'\sigma$ for some $\nu'\in\text{In}_1(\text{Aut }C)$, $\lambda'\in\text{Mult }C$.  Assertion \textbf{1} is proved. At the same time,  Assertion \textbf{2} is actually proved.~$\square$

\section*{\textbf{Chapter III. DERIVATION SPACE\\ OF INCIDENCE COALGEBRA}}\label{chapterIII}
\addtocontents{toc}{\textbf{Chapter III. DERIVATION SPACE OF INCIDENCE\\ \mbox{}\hspace{6mm}COALGEBRA}\par}

\section{On Spaces $\text{Der }(I(X,F))$}\label{section9} 

We state the main facts about the derivation space $\text{Der}(I(X,F))$ of the incidence algebra $I(X,F)$. Like before, we denote this algebra by the letter $A$.

We formulate several familiar definitions.

Let $R$ be an algebra over some commutative ring $T$. The mapping $d\colon R\to R$ is called a \textsf{derivation} of the algebra $R$ if $d$ is a linear mapping, i.e. an endomorphism of the $T$-module $R$, and the relation $d(ab)=d(a)b+ad(b)$ holds for any $a,b\in R$. All derivations of the algebra $R$ form a $T$-module. We denote it by $\text{Der }R$.

For an element $c\in R$, we define a mapping $d_c$ from $R$ into $R$ by assuming that $d_c(a)=ac-ca$, $a\in R$. Then $d_c$ is a derivation which is said to be \textsf{inner}. We say that $d_c$ \textsf{is determined} by the element $c$. All inner derivations of the algebra $R$ form a submodule of the $T$-module $\text{Der }R$. We use the symbol $\text{In(Der }R)$ to designate it.

There is a notion of a derivation in a more general form. Let $N$ be an $R$-$R$-bimodule. A homomorphism of $T$-modules $d\colon R\to N$, such that $d(ab)=d(a)b+ad(b)$ for all $a,b\in R$, is called a \textsf{derivation} of the algebra $R$ with values in the bimodule $N$. Such a derivation is said to be \textsf{inner} if there exists an element $c\in N$ such that $d(a)=ac-ca$, $a\in R$.

We use the material and notation contained at the beginning of Section 4. Thus, the splitting extension $A=L_1\oplus M_1$  will be useful.

We take an arbitrary derivation $d$ of the algebra $A$. Similar to the case of automorphisms, the derivation $d$ can be associated with the matrix $\left(\begin{smallmatrix}\alpha&\gamma\\ \delta&\beta\end{smallmatrix}\right)$ with respect to the direct decomposition $A=L_1\oplus M_1$. It is important that $\gamma=0$ by \cite[Lemma 14.1]{KryT23}\label{4}. We have the following assertion \cite[Corollaries 14.3, 15.4]{KryT23}\label{4}.

\textbf{Corollary 9.1.}
\begin{enumerate}
\item[\textbf{1.}]
Every derivation $d$ of the incidence algebra $A$ is of the form $\left(\begin{smallmatrix}\alpha&0\\ \delta&\beta\end{smallmatrix}\right)$. Here $\alpha$ is a derivation of the algebra $L_1$, $\beta$ is a derivation of the algebra $M_1$ (as a non-unital algebra), and $\delta$ is a derivation of the algebra $L_1$ with values in $L_1$-$L_1$-bimodule $M_1$.
\item[\textbf{2.}]
If $\delta=0$, then the relations
$$
\beta(ad)=\alpha(a)d+a\beta(d),\quad \beta(cb)=\beta(c)b+c\alpha(b)
$$ 
\end{enumerate}
hold for all $a,b\in L_1$ and $c,d\in M_1$.

We denote by $\text{In}_0(\text{Der }A)$ (resp., $\text{In}_1(\text{Der }A)$) the subspace of inner derivations of the algebra $A$ determined by elements of $L_1$ (resp., $M_1$).

\textbf{Lemma 9.2 \cite[Lemma 15.1]{KryT23}\label{4}.} There is a direct decomposition of spaces
$$
\text{In(\text{Der }}A)=\text{In}_0(\text{Der }A)\oplus\text{In}_1(\text{Der }A).
$$

Let the symbol $\text{Add }A$ denote the subspace in $\text{Der }A$, consisting of derivations of the form $\left(\begin{smallmatrix}0&0\\ 0&\beta\end{smallmatrix}\right)$. Such derivations are said to be \textsf{additive}. Similar to multiplicative automorphisms, we can also define them with the sue of certain systems of elements of the field $F$ (see \cite{KryT23}\label{4}).

\textbf{Definition 9.3.} Let we have a system of elements $\{c_{xy}\in F\,|\,x<y\}$ such that the relation
$$
c_{xy}=c_{xz}+c_{zy} \eqno(1)
$$
holds for all $x,z,y\in X$ with $x<z<y$. We call such systems \textsf{additive systems}. 

We can also define an additive system as a system $\{c_{xy}\in F\,|\,x\le y\}$ such that $c_{xx}=0$ for all $x$ and $c_{xy}=c_{xz}+c_{zy}$ if $x\le z\le y$. 

If $d=\left(\begin{smallmatrix}0&0\\ 0&\beta\end{smallmatrix}\right)\in\text{Add }A$, then for any $x,y\in X$ with $x<y$, there exists an element $c_{xy}\in F$ such that $\beta(b)=c_{xy}b$, where $b$ is an arbitrary element of $M_{xy}$. In this case, the system of the elements $\{c_{xy}\,|\,x< y\}$ is additive.

It turns out that this additive derivation $d$ can be associated with an  additive system of elements $c_{xy}$ ($x,y\in X$, $x<y$) of the field $F$. Conversely, every additive system of elements $\{c_{xy}\in F\,|\,x<y\}$ leads to an additive derivation. More precisely if $g=(g_{xy})\in M_1$, then we have to set
$$
d(g)=(c_{xy}g_{xy}) \text{ and } d(f)=0 \text{ for } f\in L_1.
$$

Based on the above, we write down the following statement.

\textbf{Proposition 9.4.} There is a one-to-one correspondence between additive derivations of the algebra $A$ and additive systems elements of the field $F$.

At the end of the section, we write down a theorem which describes the structure of the space $\text{Der }A$.

\textbf{Theorem 9.5 \cite{SpiO97}\label{7}.} We have the relation
$$
\text{Der }A=\text{In}_1(\text{Der }A)\oplus\text{Add }A.
$$

\section{Derivation Space $\text{Der }C$}\label{section10} 

For a while let $C$ be an arbitrary coalgebra $(C,\Delta,\varepsilon)$.

\textbf{Definition 10.1.} A linear mapping $d\colon C\to C$ is called a \textsf{derivation} of the coalgebra $C$ if the relation $\Delta d=(d\otimes 1)\Delta +(1\otimes d)\Delta $ holds.

All derivations of the coalgebra $C$ form an $F$-space which is called the \textsf{derivation space} of the coalgebra $C$. We fix  the symbol $\text{Der }C$ to denote it.

\textbf{Lemma 10.2.} If $d\in\text{Der }C$, then $d^*\in\text{Der }C^*$.

\textbf{Proof.} We have to verify that the relation
$$
d^*m=m(d^*\otimes 1) +m(1\otimes d^*)
$$
is true (mappings $m$ and $\rho$, found below, are defined at the beginning of Section 2).

It follows from the relation $\Delta d=(d\otimes 1)\Delta +(1\otimes d)\Delta $ that
$$
d^*\Delta^*=\Delta^*(d\otimes 1)^*+\Delta^*(1\otimes d)^*.
$$
We also have the relations
$$
d^*m=d^*\Delta^*\rho=\Delta^*(d\otimes 1)^*\rho+\Delta^*(1\otimes d)^*\rho,
$$
$$
m(d^*\otimes 1)+m(1\otimes d^*)=\Delta^*\rho(d^*\otimes 1)+
\Delta^*\rho(1\otimes d^*).
$$

Calculations show that the relations
$$
(d\otimes 1)^*\rho=\rho(d^*\otimes 1),\quad (1\otimes d)^*\rho=\rho(1\otimes d^*) 
$$
hold; this implies the required result.~$\square$

\textbf{Beginning with this place and until the end of the chapter, $C$ denotes the incidence coalgebra $\text{Co}(X,F)$.}

It follows from Lemma 10.2 that there are group monomorphisms
$$
\Gamma\colon \text{Der }C\to\text{Der }C^*,\quad d\to d^*,
$$
$$
\Theta\colon \text{Der }C\to\text{Der }A,
$$
where $A$ is the incidence algebra $I(X,F)$. Here we also need to take into account that the isomorphism $\Phi^*$ maps $\text{Der }C^*$ onto $\text{Der }A$. (The mapping $\Gamma$ appeared at the end of Section 2, and $\Phi^*$ and $\Theta$ appeared at the beginning of Section 3.) Theorem 12.1 actually states that in reality $\Gamma$ and $\Theta$ are isomorphisms.

In the previous section, we defined additive derivations and additive systems of elements corresponding to them (see Proposition 9.4). Now we will carry out similar considerations regarding the coalgebra $C$.

Let we have an additive system of elements $\{c_{xy}\in F\,|\,x<y\}$ (see Definition 9.3). We define a mapping $\lambda\colon C\to C$ by setting $\lambda([x,y])=c_{xy}[x,y]$ for every basis vector $[x,y]$ of the space $C$ with $x<y$, and $\lambda([x,x])=0$ for every $x$. Then $\lambda$ is a derivation of the coalgebra $C$. We call it the \textsf{additive derivation} corresponding to the additive system $\{c_{xy}\in F\,|\,x<y\}$.

All additive derivations form a space which is denoted by $\text{Add }C$.

\textbf{Proposition 10.3.} The spaces $\text{Add }A$ and $\text{Add }C$ are isomorphic. The mapping $\Theta$ is an isomorphism of these spaces.

\textbf{Proof.} We can carry out reasoning similar to the reasoning from the proof of Proposition 6.1.~$\square$

\section{Inner Derivations of Coalgebra $\text{Co}(X,F)$}\label{section11} 

We introduce the notion of an inner derivation of the incidence coalgebra $C=\text{Co}(X,F)$. The following definition also makes sense for an arbitrary coalgebra.

\textbf{Definition 11.1.} A derivation $\nu$ of the coalgebra $C$ is said to be \textsf{inner} if $\nu^*$ is an inner derivation the dual of the algebra $C^*$.

All inner derivations of the coalgebra $C$ form a subspace; we denote it by $\text{In(Der }C)$.

Let $\nu\in\text{In(Der }C)$. It is easy to verify that $\Phi^*(\nu^*)\in\text{In(Der }A)$ in our case. Therefore, we have an isomorphic embedding of spaces
$$
\Theta\colon \text{In(Der }C)\to \text{In(Der }A).
$$
The purpose of further consideration is to show that these spaces are isomorphic. Therefore, it can be argued that a derivation $\nu$ of the coalgebra $C$ is inner if and only if $\Theta(\nu)$ is an inner derivation of the algebra $A$.

Let $g$ be some function from the ring $A$ and let $\mu$ be the inner derivation determined by this function. 

For any two elements $x,y\in X$ with $x\le y$, we set
$$
\nu([x,y])=\sum_{x\le u\le y}[x,u]g(u,y)-
\sum_{x\le v\le y}g(x,v)[v,y].\eqno(1)
$$
We obtain a linear mapping $\nu$ of the space $C$.

\textbf{Proposition 11.2.} The mapping $\nu$ is a derivation of the space $C$. In addition, $\Theta(\nu)=\mu$.

\textbf{Proof.} We have to verify that the relation $\Delta\nu=(\nu\otimes 1)\Delta+(1\otimes \nu)\Delta$ is true.

We take an arbitrary basis vector $[x,y]$ of the space $C$. By calculating, we obtain that there are relations
$$
\Delta\nu([x,y])=\sum_{x\le s\le u\le y}g(u,y)([x,s]\otimes[s,u])-\sum_{x\le v\le t\le y}g(x,v)([v,t]\otimes[t,y]),\eqno(2)
$$
$$
((\nu\otimes 1)\Delta+(1\otimes\nu)\Delta)([x,y])=
$$
$$
\left[\sum_{x\le v\le p\le y}g(p,y)([x,v]\otimes[v,p])-
\sum_{x\le \ell\le z\le y}g(x,\ell)([\ell,z]\otimes[z,y])\right]+
$$
$$
\eqno(3)
$$
$$
+\left[\sum_{x\le k\le z\le y}g(k,z)([x,k]\otimes[z,y])-
\sum_{x\le v\le q\le y}g(v,q)([x,v]\otimes[q,y])\right].
$$

The expression from first pair of square brackets in $(3)$ coincides with the right part of relation $(2)$, and the expression in the second pair of square brackets is equal to zero.

Why $\Theta(\nu)=\mu$? Let $f\in A$ and let $x,y$ be two elements of $X$ with $x\le y$. With the use of Proposition 3.1, we obtain from $(1)$ the relation
$$
((\Theta(\nu))(f))(x,y)=\sum_{x\le u\le y}g(u,y)f(x,u)-
\sum_{x\le v\le y}g(x,v)f(v,y).
$$
Its right part coincides with the right part of the relation $(\mu(f))(x,y)=(fg-gf)(x,y)$.~$\square$

By the use of Proposition 10.3 and 11.2, we can write down the following result.

\textbf{Corollary 11.3.} There exists a group isomorphism
$$
\Theta\colon \text{In(Der }C)\to \text{In(Der }A).
$$

We obtain that every inner derivation of the coalgebra $C$ is of the form specified in $(1)$.

The subgroups of inner automorphisms $\text{In}_0(\text{Aut }C)$ and $\text{In}_1(\text{Aut }C)$ were defined before Corollary 7.5. We can similarly define the subspaces of inner derivations $\text{In}_0(\text{Der }C)$ and $\text{In}_1(\text{Der }C)$ of the coalgebra $C$. Then the following assertion follows from Lemma 9.2 and Corollary 11.3.

\textbf{Corollary 11.4.} We have a direct decomposition of $F$-spaces:
$$
\text{In(Der }C)=\text{In}_0(\text{Der }C)\oplus\text{In}_1(\text{Der }C).
$$

\section{Description of Derivation Space $\text{Der }C$}\label{section12} 

We write down a theorem which contain complete information about the structure of the derivation space of the coalgebra $C$.

\textbf{Theorem 12.1}. Let $A$ be an incidence algebra $I(X,F)$ and let $C$ be the incidence coalgebra $\text{Co}(X,F)$.
\begin{enumerate}
\item[\textbf{1.}]
There are relations
$$
\text{Der }C=\text{In(Der }C)+\text{Add }C=\text{In}_1(\text{Der }C)\oplus\text{Add }C.
$$
\item[\textbf{2.}]
The derivation spaces $\text{Der }C$ and $\text{Der }A$ are isomorphic.
\end{enumerate}

\textbf{Proof.} \textbf{1.} We take an arbitrary derivation $d$ of the coalgebra $C$. Then $\Theta(d)\in \text{Der }A$; see Lemma 10.2 and the text after its proof. We can write down $\Theta(d)=\mu+\psi$, where $\mu\in \text{In}_1(\text{Der }A)$, $\psi\in\text{Add }A$ (Theorem 9.5). In view of Corollary 11.3 and Proposition 10.3, there are derivations $\nu\in\text{In}_1(\text{Der }C)$ and $\lambda\in\text{Add }C$ such that $\Theta(\nu)=\mu$ and $\Theta(\lambda)=\psi$. Therefore, $\Theta(\nu +\lambda)=\mu+\psi=\Theta(d)$ and, therefore, $d=\nu +\lambda$.

\textbf{2.} The fact that $\Theta$ is the required isomorphism  follows from the proof of assertion \textbf{1} (or from Proposition 10.3 and Corollary 11.3).~$\square$

\textbf{Open Questions.}

\textbf{1.} Is it possible to define an inner automorphism of an arbitrary coalgebra $C$ (in particular, of the incidence coalgebra $C$) in terms of the coalgebra itself, i.e. without referring to the dual algebra $C^*$?

\textbf{2.} A similar question makes sense regarding inner derivations of the coalgebra $C$ (see Definitions 7.1 and 11.1).

\label{biblio}

\end{document}